\documentclass[preprint,12pt]{elsarticle}

\usepackage{psfrag}
\usepackage{graphicx,amsmath}
\usepackage{color,amsfonts,latexsym}
\usepackage{bbold}

\usepackage{eucal}
\usepackage{booktabs}
\usepackage{textcomp}
\usepackage{soul}
\usepackage{upgreek}

\usepackage{graphicx}
\usepackage{epsfig,graphics,subfigure,psfrag,amsmath,amssymb}
\usepackage{tikz}
\usepackage{multirow}
\usepackage{graphics}
\usepackage{epsfig}
\usepackage{mathptmx}
\usepackage{times}
\usepackage{amsmath}
\usepackage{amssymb}
\usepackage{indentfirst}
\newtheorem{theorem}{Theorem}
\newtheorem{lemma}{Lemma}
\newtheorem{remark}{Remark}

\newcommand{\vect}[1]{\mathrm{vec}\left(#1\right)}

\newcommand{\tr}{\mathsf{T}}

\journal{Journal of Process Control}
\date{DRAFT}

\begin{document}

\begin{frontmatter}

\title{Real-time Algorithm for Self-Reflective Model Predictive Control}

\author{Xuhui Feng, Boris Houska}

\address{School of Information Science and Technology,
ShanghaiTech University, China. \\
{\tt borish@shanghaitech.edu.cn}
}

\begin{abstract}
This paper is about a real-time model predictive control (MPC) algorithm for a particular class of model based controllers, whose objective consists of a nominal tracking objective and an additional learning objective. Here, the construction of the learning term is based on economic optimal experiment design criteria. It is added to the MPC objective in order to excite the system from time-to-time on purpose in order to improve the accuracy of the state and parameter estimates in the presence of incomplete or noise affected measurements. A particular focus of this paper is on so-called self-reflective model predictive control schemes, which have the property that the additional learning term can be interpreted as the expected loss of optimality of the controller in the presence of random measurement errors. The main contribution of this paper is a formulation-tailored algorithm, which exploits the particular structure of self-reflective MPC problems in order to speed-up the online computation. It is shown that, in contrast to generic state-of-the-art optimal control problem solvers, the proposed algorithm can solve the self-reflective optimization problems with reasonable additional computational effort and in real-time. The advantages of the proposed real-time scheme are illustrated by applying the algorithm to a nonlinear process control problem in the presence of measurement errors and process noise. 
\end{abstract}


\end{frontmatter}

\section{Introduction}
The standard variant of model predictive control (MPC)~\cite{Mayne2000,Rawlings2009} relies on the principle of \textit{certainty equivalence}: at every sampling time a nominal control performance objective is optimized subject to dynamic model equations as well as control and state constraints over a finite prediction horizon under the assumption that there are neither state estimation errors, nor external disturbances, nor any kind of model plant mismatches present, although all these errors and disturbances are the reason why a feedback controller is needed in the first place. After the MPC controller sends its first control input to the real process, the next optimization problem is solved by using the latest state estimate in order to close the loop. The success of such certainty equivalent model predictive controllers, also in industrial applications~\cite{Qin2003}, is to a large part due to the availability of fast and reliable real-time optimal control problem solvers~\cite{Diehl2002,Zavala2009}. During the last years algorithms as well as mature automatic code generation based software have been developed, which can solve nonlinear model predictive control problems online and within sampling times in the milli- and microsecond range~\cite{Houska2011,Mattingley2009}.

\bigskip
\noindent
Although one might argue that standard MPC and its more traditional variants do not attempt to achieve a tradeoff between nominal performance as well as learning objectives, many other controllers, which implement such tradeoffs, have a longer history. In fact, during the last 50 years, there have been many suggestions on how to develop controllers that implement such a tradeoff between control performance and learning. One of the pioneers of the so-called dual control problem is A.~Feldbaum, who published a whole series of seminal papers on this topic~\cite{Feldbaum1961}. Feldbaum's original dual control problem formulation is based on an optimization problem that minimizes the expected control performance, typically a least-squares tracking term, under the assumption that the accuracy of the parameter estimates depends on the information content of future state measurements. Unfortunately, explicit solutions for this original problem formulation of the dual control problem have so far only been found for very simple problems with one state variable. Numerical algorithms for solving the dual control problem in higher dimensional spaces, e.g., based on approximate dynamic programming, turn out to be rather expensive~\cite{Lee2009}. For an overview about other attempts to solve the dual control problem approximately by using techniques from the field of adaptive control the reader is referred to the overview articles~\cite{Filatov2004,Wittenmark1995}.

\bigskip
\noindent
As mentioned above, earlier or more traditional variants of MPC usually do not analyze learning objectives explicitly. However, during the last decade this situation has changed and, especially in recent years, there have appeared a significant number of articles about MPC variants that incorporate additional learning terms in order to achieve better future state and parameter estimates. For example, in~\cite{Hovd2005} it is suggested to augment the standard MPC objective by an additional term that penalizes an approximation of the variance of future state estimates. This can be implemented by augmenting the model equations with an extended Kalman filter that can be used to predict the variance of future state and parameter estimates in a linear approximation. Similar extensions of MPC with learning terms have been proposed in~\cite{Heirung2012,Heirung2015}, which augment the nominal MPC objective with optimal experiment design objectives that penalize the predicted variance of the system parameters. In recent years there have appeared a number of articles on persistently exciting MPC~\cite{Hernandez2015,Marafioti2014,Mesbah2015,Shouche2002,Zacekova2013}, which all discuss different ways to excite model predictive controllers in order to improve the accuracy of future state and parameter estimates. An even more recent trend is to extend the concept of  application oriented optimal experiment design~\cite{Hjalmarsson2009} by associated terms in the objective or constraints of an MPC problem~\cite{Larrson2013,Larsson2015}. Moreover, a recent paper on self-reflective model predictive control~\cite{Houska2016} proposes a model predictive controller that minimizes a prediction of its own expected loss of control performance in the presence of measurement noise. Notice that all these developments are closely related to the research on output-feedback MPC~\cite{Findeisen2003}. For example, in \cite{Mayne2006,Mayne2009} it is suggested to augment the dynamic system with a differential equation that implements a state estimator and then use methods from the field of robust MPC to robustly control the coupled system of the controller and estimator. This leads on the one hand to a rigorous framework for bounding or approximating the influence of the state estimator in MPC but, on the other hand, using robust tube based MPC methods also adds another layer of complexity to a potentially nonlinear dynamics that is already challenging to solve.

\bigskip
\noindent
A rather apparent drawback of all the above reviewed model predictive control schemes with additional learning objectives is that they are based on introducing additional, typically matrix-valued hyperstates, which are needed for predicting the accuracy of future state and parameter estimates. This leads to an optimal control problem that is from a numerical computation perspective much more difficult to solve than the corresponding optimal control problem without learning terms. Unfortunately, real-time algorithms, which exploit the structure of the model predictive control problems with additional learning terms, are, as far as the authors are aware, not available to date---let alone implementations and software for solving these problems reliably. Therefore, a principal goal of this paper is to develop a real-time algorithm that can exploit the structure of such problems. Here, we focus on the self-reflective model predictive control formulation, which has been proposed in~\cite{Houska2016} and which is based on augmenting a nominal MPC objective with an economic optimal experiment design criterion~\cite{Houska2015}. However, the methods in this paper can also be extended for most of the other integrated experiment design MPC methods that have been reviewed above. The main contribution of this paper is the development of a novel real-time algorithm for model predictive control with additional learning objectives. In order to avoid confusion, we mention here that the formulation of such augmented MPC problems is not an original contribution of this paper, since there are many articles about how to formulate such problems as reviewed above. The main motivation of this paper is based on the fact that the practical intention of adding a (small) correction term for the purpose of learning is, at least in most applications that are known to the authors, to obtain a refinement of an existing model based control scheme. A major change of the properties of a standard MPC scheme is usually not intended when such learning refinements are added. However, at the current status of research adding such a small correction term for the purpose of learning means---at least based on the authors' numerical experience---an increase of computation time in the order of a factor of $100-1000$ compared to standard MPC, if one of the above reviewed augmented MPC formulations is implemented by using a generic optimal control problem solver. From a practical perspective such an increase of run-time is hardly acceptable and might in fact be the main reason why MPC controllers with additional learning objectives, as developed during the last decade, are not yet used more widely in the process control industry. In order to remedy this situation, the current paper develops a tailored algorithm that can solve the self-reflective MPC problem within a sampling time that amounts to less than four times the sampling time of an associated certainty-equivalent real-time MPC scheme~\cite{Diehl2002}. The corresponding implementation is based on an extension of the automatic code generation tools that are implemented by using CASADI~\cite{Andersson2010,Andersson2012} and ACADO toolkit~\cite{Houska2011}, which can be used to implement self-reflective MPC schemes with a sampling time of less than a millisecond.

\subsection{Notation}

\noindent
We use the symbol $\mathbb S_{+}^n$ to denote the set of positive semi-definite $n \times n$ matrices. Similarly, $\mathbb S_{++}^n$ denotes the set of positive definite $n \times n$ matrices. Throughout this paper, the symbol $k \in \{ 0, \ldots, N-1 \}$ is used as a running index, which, by convention, always runs from $0$ to $N-1 \in \mathbb N$. For example, if we write a discrete-time system in the form
\[
x_{k+1} = f(x_k,u_k) \; ,
\]
then this means---if not explicitly stated otherwise---that this equation should hold for all $k \in \{ 0, \ldots, N-1 \}$.

\section{Review of real-time algorithms for nonlinear MPC}
\noindent
This section introduces notation and reviews an existing real-time algorithm for standard, certainty equivalent MPC, which is the basis for the developments in this paper. 

\subsection{Certainty equivalent MPC}
\noindent
Standard nonlinear MPC~\cite{Rawlings2009} solves at each sampling time an optimization problem of the form
\begin{eqnarray}
\label{eq::NLP}
\begin{array}{cl}
\underset{x,u}{\min}& \sum_{k=0}^{N-1}  l(x_k,u_k) + m(x_N) \\[0.2cm]
\mathrm{s.t.} &
\left\{
\begin{array}{rclrcl}
x_{k+1} &=& f(x_k,u_k) \, , & x_0 &=& y  \\[0.1cm]
\underline u &\leq& u_k \; \leq \; \overline u
\end{array}
\right.
\end{array}
\end{eqnarray}
based on the current (potentially inaccurate) state estimate $y$, sends the first control input $u_0$ to the real process, and waits for the next state estimate before the loop is repeated. Here, $x_k \in \mathbb R^{n_x}$ denotes the state, $u_k \in \mathbb R^{n_u}$ the control, $f : \mathbb R^{n_x} \times \mathbb R^{n_u} \rightarrow \mathbb R^{n_x}$ a three-times Lipschitz-continuously differentiable right-hand side function, which satisfies $f(0,0)=0$, and $\underline u, \overline u \in \mathbb R^{n_u}$ control bounds with $\underline u < 0 < \overline u$, i.e., such that Slater's constraint qualification is satisfied. For simplicity of presentation, we assume that the stage and terminal cost 
$$l(x,u) = \frac{1}{2} \left( x^\tr Q x + u^\tr R u \right) \quad \text{and} \quad m(x) = \frac{1}{2} x^\tr P_N x$$
are convex quadratic tracking terms with given positive semi-definite weighting matrices $Q,P_N \in \mathbb S_{+}^{n_x}$ and $R \in \mathbb S_{+}^{n_u}$. However, the considerations in this paper can be extended easily for the case that $l$ and $m$ are more general economic objective terms, as long as all functions are three-times Lipschitz-continuously differentiable.

\bigskip
\noindent
Notice that practical system are often modeled by continuous-time differential equations. However, there exist mature tools from the field of direct optimal control, which can be used to discretize continuous-time optimal control problems leading to a finite dimensional nonlinear programming problem of the above form~\cite{Biegler2007,Bock2000}. More general MPC problem formulations additionally comprise state constraints as well as terminal regions~\cite{Rawlings2009}, which are however not analyzed in this paper. The reason for this decision is that state constraints \mbox{can---at} least without further \mbox{precaution---lead} to infeasibility, if uncertainties are present. As the goal of this paper is on the one hand to analyze the influence of random measurement errors, but, on the other hand, is not about robust MPC, state constraints are not considered. Notice that a robustification with respect to state constraints can be taken into account by applying the output-feedback MPC based framework~\cite{Mayne2009}, which leads however to a much more complex problem formulation, which can be expected to lead to significantly larger computation times than the methods presented in the current paper.

\subsection{Real-time algorithm}
\label{sec::RTI}

\noindent
Real-time Gauss-Newton type algorithms for nonlinear MPC~\cite{Diehl2002} are divided into a \textit{preparation phase} and a \textit{feedback phase}. During the preparation phase, the function $f$ and its derivatives
\[
A(x_k,u_k) = \dfrac{\partial f(x_k, u_k)}{\partial x} \quad \text{and} \quad B(x_k,u_k) = \dfrac{\partial f(x_k, u_k)}{\partial u}
\]
are evaluated at the predicted state trajectory $x$ and control input $u$. Modern real-time MPC implementations are using advanced algorithmic differentiation methods for evaluating these derivatives efficiently and with high precision~\cite{Griewank2000,Houska2011}. Once the state measurement $y$ becomes available, the feedback phase starts. During this feedback phase, the quadratic programming problem
\begin{eqnarray}
\label{eq::QP}
\begin{array}{cl}
\underset{\Delta x, \Delta u}{\min}& \sum_{k=0}^{N-1}  l(x_k + \Delta x_k,u_k + \Delta u_k) + m(x_N + \Delta x_N) \\[0.2cm]
\mathrm{s.t.}&
\left\{
\begin{array}{l}
\Delta x_0 = y-x_0 \\[0.1cm]
\Delta x_{k+1} = A(x_k,u_k) \Delta x_k + B(x_k,u_k) \Delta u_k + f(x_k,u_k) - x_{k+1} \\[0.1cm]
\underline u \leq u_k + \Delta u_k \leq \overline u
\end{array}
\right.
\end{array}
\end{eqnarray}
is solved recalling that $l$ and $m$ are assumed to be convex quadratic forms. Next, the updated control
$$u_\text{process} = u_0 + \Delta u_0$$
is sent out to the process. The remaining states and controls are shifted in time and updated as follows:
\begin{eqnarray}
\begin{array}{rclcrcl}
\forall k \in \{ 0, \ldots, N-1 \}, \; \; x_{k} &\leftarrow& x_{k+1} + \Delta x_{k+1} \, , & & x_N &\leftarrow& x_N + \Delta x_N \, , \notag \\[0.1cm]
\forall k \in \{ 0, \ldots, N-2 \}, \; \; u_{k} &\leftarrow& u_{k+1} + \Delta u_{k+1} \, , & \text{and}& u_{N-1} &\leftarrow& u_{N-1} + \Delta u_{N-1} \, . \notag
\end{array}
\end{eqnarray}
Notice that there are many variants of the above scheme. For example, in~\cite{Diehl2002} it is suggested to prepare the linear algebra operations that are needed for solving the QP by applying a condensing based approach that eliminates the optimization variable $\Delta x$ explicitly. This leads to a smaller but more dense QP to be solved in the feedback phase. A theoretical analysis as well as stability proofs of the above real-time algorithm for nonlinear MPC under mild additional assumptions can be found in~\cite{Diehl2002}.

\subsection{Interior-point versus active set methods}
\label{sec::IP}

\noindent
The above real-time algorithm is often applied in combination with active set methods, which can be used for solving the QP~\eqref{eq::QP} efficiently. For example, in~\cite{Ferreau2008,Houska2011} online variants of active set methods for linear and nonlinear MPC are proposed, which use tailored hot-start methods in order to further speed-up the computation time during the feedback step. However, other nonlinear MPC implementations, such as~\cite{Zavala2009}, use an interior point method based framework. Here, the main idea is to introduce the modified state cost
\[
l_{\tau}(x,u) = l(x,u) - \tau \log \left[ \left( u - \underline u \right) \left( \overline u - u \right) \right] \; ,
\]
where $\tau > 0$ is a barrier parameter. Now, problem~\eqref{eq::NLP} can be approximated by an equality-constrained NLP of the form
\begin{eqnarray}
\label{eq::NLP2}
\begin{array}{cl}
\underset{x,u}{\min}& \sum_{k=0}^{N-1}  l_\tau(x_k,u_k) + m(x_N) \\[0.2cm]
\mathrm{s.t.} &
\left\{
\begin{array}{rclrcl}
x_{k+1} &=& f(x_k,u_k) \, , & x_0 &=& y \; ,
\end{array}
\right.
\end{array}
\end{eqnarray}
which is equivalent to problem~\eqref{eq::NLP} for the limit $\tau \to 0^+$ (see~\cite{Boyd2004} for a proof). The numerical advantage of this NLP is that it does not contain inequalities. A corresponding real-time scheme proceeds analogous to Section~\ref{sec::RTI}, but the difference is that the stage cost $l_{\tau}$ is at every sampling time replaced by a second Taylor order expansion at the previous iterate in order to obtain a real-time Newton-type method that solves equality-constrained QPs in the feedback phase. As the band-structure of this QP can be exploited easily, this leads to an algorithm with complexity $\mathbf{O}(N)$, which is fortunate if $N$ is large~\cite{Zavala2009}. However, a drawback of the interior-point method based real-time scheme compared to active set methods is that hot-starts are more difficult to implement. This is due to the fact that  the parameter $\tau$ has to be adjusted during the iterations. Details about how to implement an associated line-search and how to adjust the tuning parameter $\tau$ during the iterations can be found in~\cite{Boyd2004,Wright1997,Zavala2009}. Finally, notice that both the active-set as well as the interior-point method based real-time algorithm have in common that they do not solve problem~\eqref{eq::NLP} to optimality, but merely apply one Gauss-Newton or interior point iteration per sampling time.

\section{MPC with additional learning objectives}

\noindent
In practical applications the state of a system can often not be measured directly and has to be estimated from measurements. Throughout this paper, the measurements are assumed to satisfy the equation
\begin{eqnarray}
\eta_k = C x_k + v_k \; ,
\end{eqnarray}
where $C \in \mathbb R^{n_{\eta}} \times \mathbb R^{n_x}$ is a given matrix, $v_k \in \mathbb{R}^{n_\eta}$ the random measurement error, and $\eta_k \in \mathbb{R}^{n_\eta}$ the actual measurement. In a more general situation, the dependence of $\eta_k$ on $x_k$ could also be nonlinear or additionally depend on $u_k$, but the following considerations can be extended easily for this more general case. An additional complication is that the dynamic system itself may be affected by process noise,
\begin{eqnarray}
x_{k+1} = f(x_k,u_k) + w_k \; ,
\end{eqnarray}
i.e., the variables $w_k \in \mathbb R^{n_x}$ are random variables that cannot be predicted. In the following we assume that the first and second order moments of $w_k$ and $v_k$, denoted by $W \in S_{+}^{n_x}$ and $V \in S_{++}^{n_\eta}$, are given,
\begin{eqnarray}
\begin{array}{rclrclrclrcl}
\mathbb E \{w_k \} &=& 0 \; , & \mathbb E \{v_k \} &=& 0 \; , & \mathbb E \{ w_k w_k^\tr \} &=& W \; , & \mathbb E \{ v_k v_k^\tr \} &=& V \; .
\end{array}
\end{eqnarray}
Moreover, $v_k$ and $w_k$ are assumed to have bounded support such that $\Vert v_k \Vert \le \gamma$ as well as $\Vert w_k \Vert \le \gamma$ for a given radius $\gamma > 0$.
In the introduction a large number of extensions of MPC has been reviewed, which all attempt to include additional terms in the objective or constraints of the MPC controller in order take the accuracy of future state estimates into account when optimizing the control input. The construction of these additional terms uses ideas from the field of optimal experiment design as reviewed next.

\subsection{Optimal experiment design based MPC objectives}

\noindent
Optimal experiment design is a rather mature technology and we refer to the text book~\cite{Pukelsheim1993} for a general overview. As this paper is about dynamic system, it is focussed on an optimal experiment design problem formulation that has been analyzed in~\cite{Telen2013}. Here, the main idea is to introduce an extended Kalman filter in order to predict the variance matrix of future state estimates, which is based on the discrete-time Riccati recursion\footnote{The matrix $C \Sigma_k C^\tr + V$ is invertible, as we assume that $V$ is positive definite.}
\begin{eqnarray}
\begin{array}{rcl}
\Sigma_{k+1} &=& F( x_k,u_k,\Sigma_k ) \\[0.16cm]
&=& A(x_k,u_k) \left[ \Sigma_k -  \Sigma_k C^\tr \left( C \Sigma_k C^\tr + V \right)^{-1} C \Sigma_k \right]A(x_k,u_k)^\tr + W \, .
\end{array} \hspace{-0.3cm}
\end{eqnarray}
This recursion is started with the variance $\Sigma_0 = \hat \Sigma = \mathbb E \{ (y-x_0^*)( y-x_0^*)^\tr \}$ of the current state estimate $y$, where $x_0^*$ denotes the true (but unknown) initial state. Notice that in our context, the derivative function $A$ is evaluated at the predicted states and control $(x_k,u_k)$, as we use the extended Kalman filter to predict the variances $\Sigma_k$ of future state estimates without knowing the future measurements $\eta_k$ yet. This is in contrast to the traditional way of applying the extended Kalman filter in the context of estimation, where the measurements $\eta_k$ are already available and can thus be used to first update the state estimate and then evaluate $A$ at the state estimate rather than at the predicted state.

\bigskip
\noindent
In the next step, a scalar experiment design criterion $\Psi: \mathbb S_{+}^{n_x} \to \mathbb R$ can be used to penalize large variances $\Sigma_k$ in the MPC objective. In this paper we focus on a particular design criterion based on the weighted A-criterion, $\Psi_k( \Sigma_k ) = \mathrm{Tr}( \alpha_k \Sigma_k )$, $\alpha_0,\ldots,\alpha_N \in \mathbb S_{+}^{n_x}$, but in principle also other criteria such as a weighted determinant or maximum eigenvalue (D- and E-criterion) can be used as scalar measures of the size of $\Sigma_k$. The following MPC problem formulation with additional learning objective has (in a very similar form) for the first time been proposed in~\cite{Hovd2005}:
\begin{eqnarray}
\label{eq::NLP_L}
\begin{array}{cl}
\underset{x,u,\Sigma}{\min}& \sum_{k=0}^{N-1} \left\{ l_\tau(x_k,u_k) + \mathrm{Tr}( \alpha_k \Sigma_k ) \right\} + m(x_N) + \mathrm{Tr}( \alpha_N \Sigma_N )  \\[0.2cm]
\mathrm{s.t.} &
\left\{
\begin{array}{rclrcl}
x_{k+1} &=& f(x_k,u_k) \, , & x_0 &=& y  \\[0.1cm]
\Sigma_{k+1} &=& F(x_k,u_k, \Sigma_k ) \, , & \Sigma_0 &=& \hat \Sigma \; .
\end{array}
\right.
\end{array}
\end{eqnarray}
In the above MPC formulation, the matrix-valued weights $\alpha_0,\ldots,\alpha_N \in \mathbb S_{+}^{n_x}$ can be used to tradeoff between the optimal experiment design objective and the nominal control objective: roughly, small values for the $\alpha_k$s lead to better nominal tracking performance while larger values typically lead to more excitation for the purpose of learning. However, tuning these weighting matrices by hand may be cumbersome and is, to a certain extent, ambiguous. An alternative is to automatically compute more natural tradeoff weights by using the self-reflective model predictive control scheme~\cite{Houska2016} that is reviewed next.

\subsection{Self-reflective model predictive control}

\noindent
As measurement errors and process noise cannot be predicted, they cause an inevitable loss of control performance when compared to an utopia feedforward controller, which can predict all future measurement errors and process noise. Self-reflective MPC~\cite{Houska2016} is a controller that minimizes the sum of its nominal performance and a second order approximation of its own expected future loss of optimality compared to an utopia feedback controller that can predict everything. In order to outline the corresponding online optimization problem formulation, the functions
\begin{eqnarray}
\begin{array}{rcl}
X(x,u,\Omega) &=& B(x,u)^\tr P A(x,u) + \nabla_{ux}^2 \left[l_\tau(x,u) + \lambda^\tr f(x,u) \right] \\[0.16cm]
Y(x,u,\Omega) &=& B(x,u)^\tr P B(x,u) + \nabla_{uu}^2 \left[l_\tau(x,u) + \lambda^\tr f(x,u) \right] \\[0.16cm]
Z(x,u,\Omega) &=& A(x,u)^\tr P A(x,u) + \nabla_{xx}^2 \left[l_\tau(x,u) + \lambda^\tr f(x,u) \right] \\[0.16cm]
\Phi(x,u,\Omega) &=& X(x,u,\Omega)^\tr Y(x,u,\Omega)^{-1} X(x,u,\Omega)
\end{array}
\end{eqnarray}
as well as
\begin{eqnarray}
G(x,u,\Omega) = \left(
\begin{array}{l}
\lambda^\tr A(x,u) + x^\tr Q \\[0.1cm]
Z(x,u,\Omega) - \Phi(x,u,\Omega)
\end{array}
\right)^\tr \quad \text{and} \quad
M(x) = \left(
\begin{array}{c}
x^\tr P_N \\
P_N
\end{array}
\right)^\tr \notag
\end{eqnarray}
are introduced. Here, the functions $X,Y,Z,G,$ and $M$ are defined for all $x \in \mathbb R^{n_x}$ all $u \in \mathbb R^{n_u}$, and all $\Omega = [\lambda,P] \in \mathbb R \times \mathbb S_+^{n_x}$. Notice that the evaluation of the functions $X,Y,$ and $Z$ requires the computation of second order derivatives of the Hamiltonian function $\mathcal H = \l_\tau + \lambda^\tr f$. Although there exist efficient second order algorithmic differentiation schemes~\cite{Quirynen2014}, the computation of these derivatives is expensive. Also notice that the function $l_\tau$ is strictly convex for any $\tau>0$. Thus, it follows from the second order optimality conditions for problem~\eqref{eq::NLP}, in this case a discrete-time variant of Pontryagin's maximum principle~\cite{Pontryagin1962}, that the second order derivative of the Hamiltonian $\nabla_{uu}^2 \mathcal H$ is a positive definite function in the neighborhood of a solution of~\eqref{eq::NLP}. Consequently, $Y(x,u,\Omega)$ is a positive definite (and thus invertible) matrix function in this neighborhood. Now, self-reflective MPC solves at each sampling time an optimization problem of the form
\begin{eqnarray}
\label{eq::NLP_SR}
\begin{array}{cl}
\underset{x,u,\Sigma,\Omega}{\min}& \sum_{k=0}^{N-1} \left\{ l_{\tau}(x_k,u_k) + \frac{1}{2}\mathrm{Tr}\left( \Phi(x_k,u_k,\Omega_{k+1}) \Sigma_k \right) \right\} + m(x_N) \\[0.2cm]
\mathrm{s.t.} &
\left\{
\begin{array}{rclrcl}
x_{k+1} &=& f(x_k,u_k) \, , & x_0 &=& y \\[0.1cm]
\Sigma_{k+1} &=& F(x_k,u_k, \Sigma_k ) \, ,& \Sigma_0 &=& \hat \Sigma \\[0.1cm]
\Omega_k &=& G(x_k,u_k,\Omega_{k+1})  \, ,& \Omega_N &=& M(x_N) \; .
\end{array}
\right.
\end{array}
\end{eqnarray}
At this point, one might argue that the optimization problem~\eqref{eq::NLP_SR} is very similar to problem~\eqref{eq::NLP_L} in the sense that if we set $\alpha_k = \frac{1}{2}\Phi(x_k,u_k,\Omega_{k+1})$ and $\alpha_N = 0$ in problem~\eqref{eq::NLP_L} the objective functions of the two problems coincide. However, notice that in the self-reflective problem~\eqref{eq::NLP_SR} the weights $\alpha_k = \frac{1}{2}\Phi(x_k,u_k,\Omega_{k+1})$ are non-trivial functions of $x_k$ and $u_k$ and require the computation of the ancillary backward recursion for the variables $\Omega_k$, i.e., the tradeoff between the nominal tracking performance and the additive learning term is optimized. Notice that the term $\frac{1}{2}\mathrm{Tr}\left( \Phi(x_k,u_k,\Omega_{k+1}) \Sigma_k \right)$ approximates the self-reflective MPC controller's own inherent expected loss of optimality in the presence of measurement errors and process noise up to terms of order $\mathbf{O}\left(\gamma^3 \right)$, as established in~\cite{Houska2016}; see also~\cite{Stengel1994} for a dicussion and interpretation of this term in the context of unconstrained linear stochastic control. Other advantages of the self-reflective MPC formulation as well as general economic optimal experiment design criteria are discussed in~\cite{Houska2015,Houska2016}.

\subsection{Numerical challenges}

\noindent
Unfortunately, solving the integrated experiment design based MPC problem~\eqref{eq::NLP_L} is much more expensive than solving the nominal MPC formulation~\eqref{eq::NLP2}. One reason for this increase in complexity is that problem~\eqref{eq::NLP_L} comprises
$$n_x + \frac{1}{2} n_x(n_x+1)$$
state variables assuming that the symmetry of the matrix-valued state $\Sigma_k$ is exploited. The complexity of solving the self-reflective optimization problem~\eqref{eq::NLP_SR} is even worse, as this problem comprises
$$2 n_x + n_x(n_x+1)$$
states assuming again that symmetry is already exploited. Another reason is that the forward propagation of the variance matrices $\Sigma_k$ involves the evaluation of the derivative functions $A$ and $B$ of $f$, which is often much more expensive than a nominal evaluation of $f$. The backward propagation of the state $\Omega_k$ in the self-reflective MPC problem~\eqref{eq::NLP_SR} requires second derivatives of $f$, which is expensive, too. Moreover, the stage costs of both~\eqref{eq::NLP_L} as well as~\eqref{eq::NLP_SR} are non-convex while the stage cost of~\eqref{eq::NLP2} consists of a convex tracking term and, optionally, an additive self-concordant and strictly convex barrier function. As if these numerical issues would not be enough, problem~\eqref{eq::NLP_SR} comes along with the additional complication that the term $\frac{1}{2}\mathrm{Tr}\left( \Phi(x_k,u_k,\Omega_{k+1}) \Sigma_k \right)$ couples the forward states with index $k$, namely $x_k$ and $\Sigma_k$, with the backward state, namely $\Omega_{k+1}$, with index $k+1$. This destroys the separability of the objective function and renders existing real-time MPC algorithms not directly applicable, as these existing methods rely on the separability of the objective function. From a process control engineering perspective one might argue that the main motivation for adding a learning term in the MPC objective is (if this is intended at all) that small corrections, usually in the form of minor excitations from the nominal path, are induced in order to improve future state estimates. However, if adding this minor correction comes along with major numerical difficulties and a huge increase in terms of computation time, its usefulness must be assessed critically. Therefore, the goal of this paper is to develop a tailored real-time algorithm based on the self-reflective MPC problem~\eqref{eq::NLP_SR}, which addresses the above mentioned numerical challenges in such a way that only a moderate increase in computation time is encountered when adding the learning term to the MPC objective. The hope is that this will help to promote a wider use of integrated experiment design based MPC formulations in practical and industrial process control applications.

\section{Real-time algorithm for self-reflective MPC}
\label{sec::rt-alg}

\noindent
This section develops a real-time algorithm for self-reflective MPC based on the optimization problem~\eqref{eq::NLP_SR}.

\subsection{Motivation}

\noindent
The following parametric optimization problem can be interpreted as a linearly perturbed version of the original NLP~\eqref{eq::NLP2}:
\begin{eqnarray}
\label{eq::NLP3}
\begin{array}{cl}
\underset{x,u}{\min}& \sum_{k=0}^{N-1} \left\{  l_\tau(x_k,u_k) + \sigma_k^\tr u_k \right\} + m(x_N) \\[0.2cm]
\mathrm{s.t.} &
\left\{
\begin{array}{rclrcl}
x_{k+1} &=& f(x_k,u_k) \, , & x_0 &=& y \; .
\end{array}
\right.
\end{array}
\end{eqnarray}
Here, the parameter vector $\sigma \in \mathbb R^{(N-1) n_u}$ can be used to scale the perturbation. Clearly, for a given $\sigma$ the cost for solving optimization problem~\eqref{eq::NLP3} can be expected to be basically the same as the cost for solving the NLP~\eqref{eq::NLP2}, as an additional linear term in the objective hardly adds numerical difficulties. In order to avoid confusion at this point, it is mentioned that the method presented here should not be mixed up with the so-called \textit{modifier adaption method}~\cite{Marchetti2009}, although one might argue that the parameter $\sigma$ could be interpreted as a ``modifier'' that corrects the stage cost. However, in the context of the modifier adaption method, the correction of the objective is added in order to correct model-plant mismatches and $\sigma$ is updated based on measurements of the objective gradient and constraints violation. In contrast to this approach, the intention of the presented framework is different, as we intend to perturb the gradient of the stage cost in order to improve the future state estimates while unstructured model-plant mismatches are beyond the scope of the current paper. In order to establish a connection between problem~\eqref{eq::NLP3} and the self-reflective MPC controller~\eqref{eq::NLP_SR} the following auxiliary function is introduced:
\begin{eqnarray}
\begin{array}{rccl}
E(u,x_0,\Sigma_0) &=& \underset{x,\Sigma,\Omega}{\min} & \sum_{k=0}^{N-1} \frac{1}{2}\mathrm{Tr}\left( \Phi(x_k,u_k,\Omega_{k+1}) \Sigma_k \right) \\[0.2cm]
& &\mathrm{s.t.}&
\left\{
\begin{array}{rclrcl}
x_{k+1} &=& f(x_k,u_k) \, , \\[0.1cm]
\Sigma_{k+1} &=& F(x_k,u_k, \Sigma_k ) \, , \\[0.1cm]
\Omega_k &=& G(x_k,u_k,\Omega_{k+1})  \, ,& \Omega_N &=& M(x_N) \; .
\end{array}
\right.
\end{array}
\end{eqnarray}
Next, the equivalence of the optimization problems~\eqref{eq::NLP_SR} and~\eqref{eq::NLP3} can be established under the following conditions.

\bigskip
\begin{lemma}
\label{lem::equivalence}
Let $f$ be a three-times Lipschitz-continuously differentiable function, $y = \mathbf{O}(\gamma)$ and $\Vert \hat \Sigma \Vert = \mathbf O (\gamma^2)$, and let $(x^*,u^*,\Sigma^*,\Omega^*)$ be a regular local minimizer of the self-reflective optimization problem~\eqref{eq::NLP_SR}. If $\gamma$ is sufficiently small, then the function $E$ is differentiable in a neighborhood of $u^*$ and the gradient of $E$ with respect to $u$ is Lipschitz continuous. Moreover, if we set $\sigma = \nabla_u E(u^*,y,\hat \Sigma)$, then $(x^*,u^*)$ is a regular local minimizer of the optimization problem~\eqref{eq::NLP3}.
\end{lemma}

\bigskip
\noindent
\textit{Proof.}
Let us introduce the auxiliary function
\begin{eqnarray}
\begin{array}{rcl}
L(u,x_0) \; = &\underset{x}{\min}& \sum_{k=0}^{N-1} l_\tau(x_k,u_k) + m(x_N) \\[0.2cm]
&\mathrm{s.t.} &
\left\{
\begin{array}{rclrcl}
x_{k+1} &=& f(x_k,u_k) \; .
\end{array}
\right.
\end{array}
\end{eqnarray}
such that problem~\eqref{eq::NLP_SR} can be written in the form $\min_u \left\{ L(u,y) + E(u,y,\hat \Sigma) \right\}$. The fact that $E$ is Lipschitz-differentiable with respect to $u$ follows from the fact that $f$ is assumed to be three-times Lipschitz-continuously differentiable and the assumption that $\gamma$ is sufficiently small, which implies that the function $Y$ is invertible in the neighborhood of a regular minimizer as established in~\cite{Houska2016}. Moreover, $L$ is three times continuously differentiable in $u$. Consequently, the first order necessary optimality condition
\[
0 \, = \, \nabla_u L(u^*,y) + \nabla_u E(u^*,y,\hat \Sigma)  \quad \Longrightarrow \quad \sigma \, = \, - \nabla_u L(u^*,y) \; .
\]
must be satisfied, which implies that $u^*$ is a stationary point of Problem~\eqref{eq::NLP3}. Moreover, as $l_\tau$ is a strictly convex function and $f(0,0)=0$, the Hessian matrix $\nabla_u^2 L(0,0)$ is positive definite, which in turn implies that
$$\nabla_u^2 L(u^*,y) \, = \, \nabla_u^2 L(0,0) + \mathbf{O}(\gamma)$$
is positive definite, too, for sufficiently small $\gamma$. This follows from the fact that $L$ is three times continuously differentiable. Thus, $u^*$ is a regular local minimizer of the optimization problem
\[
\min_{u} \; \left\{ L(u,y) + \sigma^\tr u \right\} \; ,
\]
which is the statement of the lemma.\qed

\bigskip
\noindent
Notice that Lemma~\ref{lem::equivalence} is only a technical intermediate result, which is in this form not useful yet, as it relies on the availability of the (in a real-time context not available) optimal solution $u^*$ of the optimization problem~\eqref{eq::NLP_SR} in order to compute the optimal perturbation vector $\sigma$. However, the main idea of the current paper is to develop a real-time algorithm that computes gradient updates of the form $\sigma = \nabla_u E(u,y,\hat \Sigma)$ at a suitable real-time control iterate $u$ rather than computing the optimal perturbation vector $\nabla_u E(u^*,y,\hat \Sigma)$. Before we analyze this idea in more detail, the following section first focusses on an algorithm for evaluating this gradient efficiently and in real-time mode.

\subsection{Fast computation of the gradient}

\noindent
In order to develop an efficient algorithm for evaluating the gradient of the function $E$, it is convenient to introduce the stacked notation
\begin{eqnarray}
\begin{array}{rclrcl}
\kappa_k &=& \left(
\begin{array}{c}
x_k \\
\vect{\Sigma_k}
\end{array}
\right) \, , & \mathcal F(\kappa_k,u_k) &=& \left(
\begin{array}{c}
f(x_k,u_k) \\
\vect{F(x_k,u_k,\Sigma_k)}
\end{array}
\right)
\\[0.5cm]
\omega_k &=& \vect{\Omega_k}, & \mathcal G(\kappa_k,\omega_{k+1},u_k) &=& \vect{G(x_k,u_k,\Omega_{k+1})} \\[0.5cm]
& & & \mathcal L( \kappa_k, \omega_{k+1}, u_k ) &=& \frac{1}{2}\mathrm{Tr}\left( \Phi(x_k,u_k,\Omega_{k+1}) \Sigma_k \right) \; ,
\end{array}
\end{eqnarray}
where the function ``$\mathrm{vec}$'' stacks all independent components of a matrix into a vector exploiting symmetry whenever possible.
By using this notation the function $E$ can be written in the more compact form
\begin{eqnarray}
\label{eq::Esummary}
\begin{array}{rccl}
E(u,\kappa_0) &=& \underset{\kappa,\omega}{\min} & \sum_{k=0}^{N-1} \mathcal L( \kappa_k, \omega_{k+1}, u_k ) \\[0.2cm]
& &\mathrm{s.t.}&
\left\{
\begin{array}{rclrcl}
\kappa_{k+1} &=& \mathcal F( \kappa_k, u_k ) \, , \\[0.1cm]
\omega_k &=& \mathcal G(\kappa_k,\omega_{k+1},u_k) \, ,& \omega_N &=& \mathcal M(\kappa_N) \; .
\end{array}
\right.
\end{array}
\end{eqnarray}
A four-sweep algorithm for computing the gradient of $E$ for a given input vector $u$ is outlined in Algorithm~1. 
\begin{figure}[htb]
\footnotesize
\begin{center}
 \begin{minipage}{1\textwidth}
\rule{1\textwidth}{0.3mm}\\
\textbf{Algorithm 1: Four-sweep algorithm for updating the perturbation vector $\sigma$.}\\
\rule{1\textwidth}{0.3mm}\\
\textbf{Input:} Control vector $u$; initial state $x_0$ and initial variance $\Sigma_0$.\\

\textbf{Four-sweeps (with $k \in \{ 0, \ldots, N-1 \}$):}\\[-0.3cm]
\begin{quote}
\begin{enumerate}
\setlength{\itemsep}{2pt}

\item \textit{Nominal forward sweep.} Set $\kappa_0 = [x_0, \text{vec}(\Sigma_0)]$ and iterate forwards
\[
\kappa_{k+1} = \mathcal F( \kappa_k, u_k ) \; .
\]

\item \textit{Nominal backward sweep.}  Set $\omega_N = \text{vec}(\Omega_N)$ and iterate backwards
\[
\omega_{k} = \mathcal G( \kappa_k, \omega_{k+1}, u_k ) \; .
\]

\item \textit{Adjoint forward sweep.} Set $a_0 = 0$ and iterate forwards
\[
a_{k+1} = \nabla_{\omega} \mathcal G( \kappa_k, \omega_{k+1}, u_k ) a_k  - \nabla_{\omega} \mathcal L( \kappa_k, \omega_{k+1}, u_k ) \; .
\]

\item \textit{Adjoint backward sweep.} Set $b_N = \nabla \mathcal M(\kappa_N) a_N$ and iterate backwards
\[
b_k = \nabla_{\kappa} \mathcal G( \kappa_k, \omega_{k+1}, u_k ) a_{k} + \nabla_{\kappa} \mathcal F( \kappa_k, u_k ) b_{k+1}  - \nabla_{\kappa} \mathcal L( \kappa_k, \omega_{k+1}, u_k ) \; .
\]

\item \textit{Final evaluation of the gradient.} Set
\[
\sigma_k = \nabla_{u} \mathcal L( \kappa_k, \omega_{k+1}, u_k ) - \nabla_{u} \mathcal G( \kappa_k, \omega_{k+1}, u_k ) a_k - \nabla_{u} \mathcal F( \kappa_k, u_k ) b_{k+1} \; .
\]

\end{enumerate}
\end{quote}
\smallskip

\textbf{Output:} The perturbation vector $\sigma = \left[ \sigma_0^\tr, \ldots, \sigma_{N-1}^\tr \right]^\tr$.\\
\removelastskip\rule{1\textwidth}{0.3mm}
\end{minipage}
\end{center}
\end{figure}
Notice that Algorithm~1 proposes a mixed forward-backward algorithmic differentiation scheme that exploits the particular structure of the function $E$. The derivation of this algorithm exploits the concept of duality in linear programming as summarized in the following theorem.

\bigskip
\begin{theorem}
If $f$ is three times continuously differentiable, then Algorithm~1 returns the gradient vector $\sigma = \nabla_u E(u,y,\Sigma_0)$.
\end{theorem}

\bigskip
\noindent
\textit{Proof.}
The directional forward derivative $\nabla_u E(u,y,\Sigma_0)^\tr \Delta u$ of the function $E$ in direction $\Delta u$ can be found by solving the following first order linear variation of problem~\eqref{eq::Esummary}:
\begin{eqnarray}
\label{eq::variation}
\begin{array}{l}
\nabla_u E(u,y,\Sigma_0)^\tr \Delta u = \\[0.2cm]
\begin{array}{cl}
\underset{\Delta \kappa, \Delta \omega}{\min} & \sum_{k=0}^{N-1} \left(
\begin{array}{c}
\nabla_{\kappa} \mathcal L( \kappa_k, \omega_{k+1}, u_k ) \\
\nabla_{\omega} \mathcal L( \kappa_k, \omega_{k+1}, u_k ) \\
\nabla_{u} \mathcal L( \kappa_k, \omega_{k+1}, u_k )
\end{array}
\right)^\tr
\left(
\begin{array}{c}
\Delta \kappa_k \\
\Delta \omega_{k+1} \\
\Delta u_k
\end{array}
\right) \\[0.2cm]
\mathrm{s.t.}&
\left\{
\begin{array}{rcl}
\Delta \kappa_0 &=& 0 \, , \\[0.1cm]
\Delta \kappa_{k+1} &=& \nabla_{\kappa} \mathcal F( \kappa_k, u_k )^\tr \Delta \kappa_k + \nabla_{u} \mathcal F( \kappa_k, u_k )^\tr \Delta u_k \, ,\\[0.1cm]
\Delta \omega_N &=& \nabla \mathcal M(\kappa_N)^\tr \Delta \kappa_{N} \; , \\[0.1cm]
\Delta \omega_k &=& \left(
\begin{array}{c}
\nabla_{\kappa} \mathcal G( \kappa_k, \omega_{k+1}, u_k ) \\
\nabla_{\omega} \mathcal G( \kappa_k, \omega_{k+1}, u_k ) \\
\nabla_{u} \mathcal G( \kappa_k, \omega_{k+1}, u_k )
\end{array}
\right)^\tr
\left(
\begin{array}{c}
\Delta \kappa_k \\
\Delta \omega_{k+1} \\
\Delta u_k
\end{array}
\right) \, .
\end{array}
\right.
\end{array}
\end{array}
\end{eqnarray}
By writing out the first order optimality conditions of the above linear programming problem, explicit expressions for the multiplies $a_0,a_1, \ldots, a_N$ as well as $b_N,b_{N-1},\ldots,b_0$ are found. They are given by the recursion in Step~3 and~4 of Algorithm~1. Moreover, since Problem~\eqref{eq::variation} is a linear programming problem, there is no duality gap, i.e., the objective values of the primal and dual objectives coincide,
\begin{eqnarray}
\begin{array}{l}
\nabla_u E(u,y,\Sigma_0)^\tr \Delta u = \\[0.16cm]
\sum_{k=1}^{N-1} \Delta u_k^\tr \left[ \nabla_{u} \mathcal L( \kappa_k, \omega_{k+1}, u_k ) - \nabla_{u} \mathcal G( \kappa_k, \omega_{k+1}, u_k ) a_k - \nabla_{u} \mathcal F( \kappa_k, u_k ) b_{k+1} \right] \; .
\end{array} \notag
\end{eqnarray}
As this equation holds for all directions $\Delta u$, a comparison of  coefficients yields $\sigma = \nabla_u E(u,y,\Sigma_0)$, as stated by the theorem.\qed

\subsection{Real-time algorithm}
\label{sec::Algorithm}

\noindent
A real-time iteration scheme for self-reflective MPC is obtained by implementing the following steps in a loop.
\begin{enumerate}
\item Compute the gradient $\sigma = \nabla_u E( u, x_0, \hat \Sigma )$ by using Algorithm~1.
\item Collect new measurements and use an extended Kalman filter to update $y$.
\item Compute a local minimizer $(x^+,u^+)$ of the optimization problem~\eqref{eq::NLP3}.
\item Send the solution $u_0^+$ to the process.
\item Shift all variables, $u_k \leftarrow u_{k+1}^+$, $u_{N-1}=u_{N-1}^+$, $x_k \leftarrow x_{k+1}^+$, and $x_N = x_N^+$. \item Update the variance matrix $\hat \Sigma \leftarrow F(x_0^+,u_0^+, \hat \Sigma)$.
\end{enumerate}
Notice that in the above scheme, the perturbation vector $\sigma$ is updated in the preparation phase in Step~1, i.e., before the actual measurement arrives. This is motivated by the fact that the computation of this gradient is usually the most expensive step and therefore done before the feedback phase, Step~2-4. Notice that variants of the above scheme might refine Step~3 and solve a quadratic approximation of \eqref{eq::NLP3} in order to further speed-up the feedback time, which leads to a scheme that is similar to the standard real-time algorithms in Sections~\ref{sec::RTI} and~\ref{sec::IP}. However, before we discuss such variants, let us first analyze, why the above real-time scheme can be expected to be contractive at all. For this aim, we introduce the auxiliary function\footnote{Notice that (if $L$ is not strictly convex in $v$) \eqref{eq::argmin} introduces a slight abuse of notation in the sense that this paper uses local optimization algorithms only. In the context of this paper, \eqref{eq::argmin} should be read as ``$\phi(u)$ is the local minimizer that is found by initializing the local solver in a neighborhood of the solution $u^*$ of \eqref{eq::NLP_SR}''. }
\begin{eqnarray}
\label{eq::argmin}
\phi(u) \; = \; \underset{v}{\mathrm{argmin}} \; \left\{ L(v,y) + \nabla_u E(u,y,\hat \Sigma)^\tr v \right\} \; .
\end{eqnarray}
Notice that Lemma~\ref{lem::equivalence} implies that the solution $u^*$ of the self-reflective optimization problem~\eqref{eq::NLP_SR} is a fix-point of the function $\phi$, $u^* = \phi(u^*)$. Moreover, the update from Step~3 of the above outlined real-time procedure can be written in the form
\[
u^+ = \phi(u) \; .
\]
Thus, the proposed algorithm can be interpreted as a real-time variant of a fixed point iteration. Now, a contractivity condition can be established by using Banach's fixed point theorem.

\bigskip
\noindent
\begin{theorem}
\label{thm::contraction}
Let all the conditions from Lemma~\ref{lem::equivalence} be satisfied and let the current iterate $u$ for the control input be in a neighborhood of $u^*$. The update from Step~3 of the proposed real-time scheme contracts linearly, i.e., we have
\[
\Vert u^+ - u^* \Vert \; \leq \; c \, \Vert u - u^* \Vert
\]
with contraction constant $c < 1$.
\end{theorem}

\bigskip
\noindent
\textit{Proof.}
Lemma~\ref{lem::equivalence} assumes that the local minimizer $u^*$ is regular. Consequently, the optimization problem~\eqref{eq::argmin} is regular for all $u$ in a sufficiently small neighborhood of $u^*$. By using a standard result from parametric nonlinear programming~\cite{Nocedal2006}, this implies that the map $\phi$ is locally Lipschitz continuous with respect to the perturbation gradient,
\[
\phi(u_1) -\phi(u_2) = \mathbf{O}\left( \left\| \nabla_u E(u_1,y,\hat \Sigma) - \nabla_u E(u_2,y,\hat \Sigma) \right\| \right)
\]
for all $u_1,u_2$ in a neighborhood of $u^*$. Next, we can use that the gradient $\nabla_u E$ is Lipschitz continuous, which implies
\[
\phi(u_1) -\phi(u_2) = \mathbf{O}\left( \gamma \right) \Vert u_1 - u_2 \Vert \; ,
\]
since $E$---and consequently also the Lipschitz constant of $\nabla_u E$---is of order $\mathbf{O}(\gamma)$. Thus, we have
\[
\Vert u^+ - u^* \Vert \; = \; \Vert \phi(u) - \phi(u^*) \Vert \; \leq \; c \, \Vert u - u^* \Vert \; ,
\]
with $c = \mathbf{O}(\gamma)$, i.e., the update contracts if $\gamma$ is sufficiently small.\qed

\bigskip
\noindent
As any other gradient based optimization method, also the performance of the above outlined real-time scheme can be improved further by implementing a pre-conditioner~\cite{Nocedal2006}. In this paper, the optimization problem~\eqref{eq::NLP_SR} is solved once offline at the steady-state. The corresponding Hessian matrix, $\nabla_u^2 E(0,0,\gamma^2 I)$, or a suitable approximation of this matrix, can be used to scale the optimization variable $u$ offline in order to improve the contraction rate of the proposed real-time scheme. As mentioned above, another improvement of the above real-time scheme refines Step~3 by solving a quadratic approximation of \eqref{eq::NLP3} at every sampling time instead of solving this problem to optimality. In this case, the evaluation of all gradients can be done in the preparation step, as explained in Section~\ref{sec::RTI}. A discussion of why contractivity of the real-time iterates is enough to ensure nominal local closed-loop stability of the corresponding MPC controller for exact state measurements and under mild regularity assumptions can be found in~\cite{Diehl2002}. However, in the context of this paper, we need to make the additional assumption that the system is locally observable such that the extended Kalman filter yields bounded variance matrices of order $\hat \Sigma = \mathbf{O}(\gamma^2)$; see also~\cite{Diehl2009,Haverbeke2011} for an overview of how to analyze the closed loop-stability of MPC-EKF cascades.

\section{Numerical case study}
\label{sec::case-sudy}

\subsection{Chemical process model}
\noindent
This paper analyzes a controllable chemical reaction, where the discrete-time right-hand side function $f(x,u) = z(h,x,u)$ is given in the form of the solution of the differential equation system 
\begin{eqnarray}
\frac{\partial z(t,x,u)}{\partial t} &=& \left(
\begin{array}{l}
-(D+k_1) z_1(t,x,u) - k_2 z_2(t,x,u) z_3(t,x,u) + u_1 \\
- D z_2(t,x,u) - k_3 z_2(t,x,u) z_3(t,x,u) + k_4 z_1(t,x,u) + u_2 \\
- D z_3(t,x,u) - k_5 z_2(t,x,u) z_3(t,x,u) + u_3
\end{array}
\right) \notag \\[0.15cm]
z(0,x,u) &=& x \; , \notag
\end{eqnarray}
which has to be evaluated numerically, e.g., by using a Runge-Kutta integrator. The given constant $h=\frac{1}{2}$ denotes the discrete-time step-size. The differential states $z_1,z_2$ and $z_3$ denote the concentrations of three substances, which react with each other. For simplicity of presentation, the corresponding reaction constants $k_1 = k_2 = k_4 = \frac{1}{2}$ and $k_3 = k_5 = \frac{1}{10}$ as well as the dilution rate $D = 0.1$ are assumed to be given. The feeding rates $u_1,u_2,$ and $u_3$ are control inputs, which can be used to adjust the inflow of the substances. We assume that only the concentration $z_1$ can be measured, $C = [1,0,0]^\tr$. The variance of the measurement error is given by $V = 0.005^2$. A summary of all model parameters can be found in Table~\ref{tab:parameters}.

\begin{table}[h]
\footnotesize
\begin{center}
\begin{tabular}{lcc}
\toprule
Name & Symbol & Value\\
\midrule
discrete-time step-size & $h$ & $0.5$ \\
MPC horizon length & $N$ & $20$ \\
lower bound on the control & $\underline u$ & $[0,-1,0]^\tr$ \\
upper bound on the control & $\overline u$ & $[\infty,\infty,\infty]$ (not  implemented) \\
initial state estimate & $\hat y$ & $[1,\;5,\;0]^\tr$\\
state reference & $x_\text{ref}$ & $[1,\;5,\;0]^\tr$\\
control reference & $u_\text{ref}$ & $[0.6,\; 0,\; 0]^\tr$\\
initial state variance & $\widehat{\Sigma}$ & $0$\\
measurement error variance & $V$ & $2.5 \cdot 10^{-5}$\\
process noise variance & $W$ & $\textrm{diag} \left( 0,\;0.64,\;0 \right)$\\
state weighting matrix & $Q$ & $\textrm{diag} \left( 1,\;1,\;1 \right)$\\
control weighting matrix & $R$ & $\textrm{diag} \left( 1,\;1, \;100 \right)$\\[0.2cm]
terminal cost weighting matrix & $P_N$ & $\textrm{diag} \left( 1,\;1,\;1 \right)$ \\
\bottomrule
\end{tabular}
\end{center}
\caption{\label{tab:parameters}Parameter values.}
\end{table}

\noindent
The stage and terminal cost functions are given by the quadratic expressions
$$l(x,u) = \frac{1}{2} \left( (x-x_\text{ref})^\tr Q (x-x_\text{ref}) + (u-u_\text{ref})^\tr R ( u-u_\text{ref}) \right)$$
and
$$m(x) = \frac{1}{2} (x-x_\text{ref})^\tr P_N (x-x_\text{ref}) \; .$$
At this point, we have introduced a small abuse of notation as in the previous sections we have set $x_\mathrm{ref} = 0$ and $u_\text{ref} = 0$. However, the above problem can easily be brought into this form by shifting all states and controls by a constant offset.

\subsection{Implementation details}

\noindent
The implementation in this paper is based on the real-time algorithm from Section~\ref{sec::Algorithm}. Here, a pre-conditioner is implemented by computing the exact Hessian of $E$ offline at the steady-state, as explained in Section~\ref{sec::Algorithm}, and we implement the refined version of Step~3, i.e., only a quadratic approximation of Problem~\eqref{eq::NLP3} is solved at every sampling time. The barrier parameter $\tau$ is kept constant and set to $\tau = 0.001$, which turns out to lead to sufficiently accurate results for this particular case study. Moreover, we use automatic code generation tools in order to export all algorithmic routines in the form of optimized C-code~\cite{Houska2011}. In particular, the four-sweep method for the perturbation vector update, Algorithm~1, is implemented by using the algorithmic differentiation software~\texttt{CASADI}~\cite{Andersson2010,Andersson2012}. Additionally, we use the software \texttt{qpOASES}~\cite{Ferreau2008} as a QP solver for implementing the real-time update.  All the results below are obtained on a Mac OX EI Captain operating system with $2.6$ GHz processor and $8$ GB, $1600$ MHz DDR3.

\begin{figure}
\begin{center}
\includegraphics[scale=0.58]{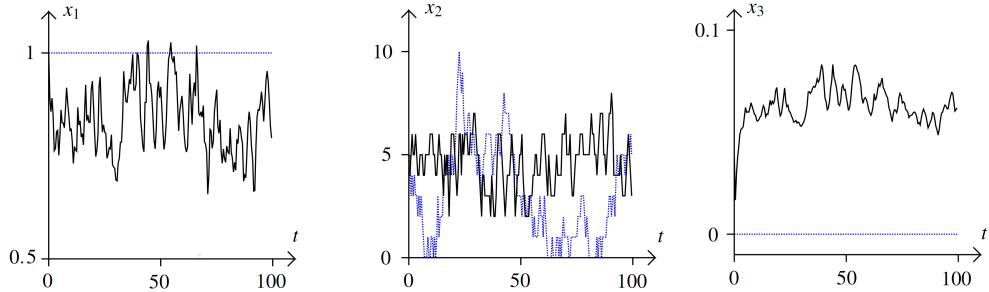}
\end{center}
\caption{\label{fig::performance} Closed-loop control trajectories of certainty-equivalent MPC (blue, dotted) and self-reflective MPC (black, solid) in the presence of random process noise and measurement errors.}
\end{figure}

\subsection{Control performance}

\noindent
Figure~\ref{fig::performance} shows a comparison of the closed loop state trajectories for randomly chosen uncertainty scenarios. The blue dotted lines in the left, middle, and right plot show the closed loop trajectories that have been obtained by running a certainty equivalent MPC controller in combination with an extended Kalman filter. We have simulated uniformly distributed process noise and measurement errors. The corresponding variance matrices $W$ and $V$ are listed in Table~\ref{tab:parameters}. Notice that the certainty equivalent MPC controller keeps the states $x_1$ and $x_3$ very close to their steady-state values, as this controller does not detect a need to excite the system. This is contrast to the self-reflective MPC controller, which excites the states $x_3$ and $x_1$ in order to be able to estimate all states more accurately. Recall that only the first state $x_1$ can be measured in our example and thus an excitation is necessary in order to be able to gather information about the other states. The main difference between the certainty equivalent MPC controller and the self-relfective MPC controller can be seen when looking at the state $x_2$. While the self-reflective MPC controller is able to keep this state in a neighborhood of its reference state $[x_\mathrm{ref}]_2 = 5$, the corresponding closed-loop trajectory of the certainty equivalent controller leads to much larger deviations in the presence of noise, as almost no information about this state is available when keeping $x_1$ and $x_3$ close to their steady-state values. This difference between the two controllers can be illustrated by analyzing their average performance. For the certainty-equivalent MPC controller the value
\[
\frac{1}{M} \, \sum_{i=1}^M l \left( x_i^{\mathrm{MPC}},u_i^{\mathrm{MPC}} \right) \, \approx \, 11.4
\]
is found for sufficiently long closed-loop simulation times $M$. The corresponding average performance of the self-reflective MPC controller is given by
\[
\frac{1}{M} \, \sum_{i=1}^M l \left( x_i^{\mathrm{SRMPC}},u_i^{\mathrm{SRMPC}} \right) \, \approx \, 2.2
\]
using the same random uncertainty scenario. Of course, the main reason why the performance difference between the two controllers is for this example particularly large is mainly due to the fact that the considered control system is not observable at its steady-state such that the self-reflective excitations lead to a major performance improvement in the presense of noise. For other applications the improvement of control performance may be less spectacular, but in general self-reflective MPC performs better than certainty-equivalent MPC and therefore, at least for some applications, this type of self-relfective controllers should be of practical interest.

\subsection{Runtime performance}

\noindent
Table~\ref{tab:runTime} summarizes the run-time of the different steps of the proposed real-time algorithm for self-reflective MPC.
\begin{table}[h]
\footnotesize
\begin{center}
\begin{tabular}{lrr}
\toprule
& CPU time ($\upmu$s) & \%\\
\midrule
Computation of the perturbation vector update (Algorithm~1) & $325$ & $68$ \\
Preparation of the QP (evaluation of $A$, $B$, etc.) & $41$ & $9$ \\
Computation of $x^+$ and $u^+$ with qpOASES (feedback step) & $108$ & $22$ \\
Other operations (including Kalman filter update) & $\leq 3$ & $\approx 1$ \\
\midrule
Total time & $477$ & $100$ \\
\midrule
Comparison: total CPU time of real-time MPC without learning term & $152$ & $32$ \\
\bottomrule
\end{tabular}
\end{center}
\caption{\label{tab:runTime} CPU time of one real-time step.}
\end{table}
The computation of the perturbation vector update takes approximately $325 \upmu$s, which corresponds to about $68 \%$ of the overall CPU time needed for one real-time step of the proposed self-reflective MPC algorithm. A standard MPC controller without learning terms needs $152 \upmu$s per real-time step, as such a standard controller needs to implement basically the same operations except for the perturbation vector update. Given the fact that solving the original self-reflective optimization problem~\eqref{eq::NLP_SR} with a generic optimal control solver takes more than $800$ times longer than the standard MPC problem, the above run-times must be considered as a major improvement compared to the state-of-the-art, although one might still argue that an increase of a factor $3$-$4$ in terms of run-time for adding a self-reflective learning term correction is still a lot. 

\begin{remark}
The case study in this paper has only three states and must be considered a small-scale example. Thus, there arises the question how the proposed real-time algorithm performs for larger scale examples. Although we do not present such larger case studies in detail as part of this paper, preliminary numerical experiments indicate that also for larger problems one real-time step of the proposed algorithm takes approximately three to four times longer than a standard real-time MPC step, i.e., also for larger problems a relatively moderate increase in run-time is observed if the additive self-reflective learning term is taken into account. Such a behavior can also be expected from our theoretical considerations in the sense that both the standard MPC real-time step applied to a problem formulation without learning term as well as Algorithm~1 have a complexity of order $\mathbf{O}( N n_x^3 )$ recalling that $n_x$ denotes the state dimension and $N$ the horizon length.\hfill$\diamond$
\end{remark}

\section{Conclusion}
\label{sec::conclusion}

\noindent
This paper has proposed a novel real-time algorithm for integrated experiment design MPC. This algorithm improves the run-time performance compared to existing generic real-time MPC algorithms by orders of magnitude, as it is capable of exploiting the particular structure of MPC problems with additional learning terms. A particular focus of this paper has been on self-reflective MPC, a controller whose objective is to minimize the sum of a nominal tracking term and the expected loss of optimality of its own performance in the presence of random process noise and measurement errors. Lemma~\ref{lem::equivalence} has established that such self-reflective model predictive controllers are equivalent to a standard MPC problem with affinely perturbed objective function as long as the perturbation vector is chosen to be equal to the gradient of the expected loss of optimality. As this affinely perturbed MPC problem can be solved as fast as standard MPC problems without learning terms, the only additional cost is the cost of computing or approximating the optimal perturbation vector. A four-sweep algorithm for computing this perturbation vector approximately and in real-time has been proposed in Algorithm~1. The properties of the associated overall real-time scheme for self-reflective MPC have been analyzed in Theorem~\ref{thm::contraction} and have been illustrated numerically by applying it to a nonlinear predator-prey-feeding control problem. For the presented case study, one real-time step takes $325 \, \upmu$s , approximately three times more than a standard real-time MPC step without learning perturbation. However, compared to existing generic real-time MPC methods applied to the self-reflective problem this amounts to a moderate increase in run-time that may be acceptable in practice, if the additional learning terms leads---as in the presented example---to a significantly improved overall control performance in the presence of measurement errors and process noise.

\section*{Acknowledgements}

\footnotesize
\noindent
This research was supported by National Natural Science Foundation China, \mbox{Nr.~61473185}, as well as ShanghaiTech University, Grant-Nr.~\mbox{F-0203-14-012}.

\normalsize
\section*{References}

\bibliographystyle{plain}

\begin{thebibliography}{99}

\bibitem{Andersson2010}
J.~Andersson, B.~Houska, M.~Diehl.
\newblock Towards a Computer Algebra System with Automatic Differentiation for Use with Object-Oriented Modelling Languages.
\newblock {\em 3rd International Workshop on Equation-Based Object-Oriented Modeling Languages and Tools}, pp.~99--105, Oslo, Norway, October~3, 2010.

\bibitem{Andersson2012}
J.~Andersson, J.~Akesson, M.~Diehl.
\newblock CasADi: A Symbolic Package for Automatic Differentiation and Optimal Control.
\newblock In {\em Recent Advances in Algorithmic Differentiation}, Volume 87 of the series Lecture Notes in Computational Science and Engineering, p.~297--307, 2012.

\bibitem{Biegler2007}
L.T. Biegler.
\newblock {A}n overview of simultaneous strategies for dynamic optimization.
\newblock {\em Chemical Engineering and Processing}, 46:1043--1053, 2007.

\bibitem{Bock2000}
H.G. Bock, M.~Diehl, D.B. Leineweber, and J.P. Schl\"oder.
\newblock {A} direct multiple shooting method for real-time optimization of
  nonlinear {DAE} processes.
\newblock In F.~Allg\"ower and A.~Zheng, editors, {\em Nonlinear Predictive
  Control}, volume~26 of {\em Progress in Systems Theory}, pages 246--267,
  Basel Boston Berlin, 2000. Birkh\"auser.

\bibitem{Boyd2004}
S.~Boyd and L.~Vandenberghe.
\newblock {\em {C}onvex {O}ptimization}.
\newblock University {P}ress, Cambridge, 2004.

\bibitem{Diehl2002}
M.~Diehl, H.G. Bock, J.P. Schl\"oder, R.~Findeisen, Z.~Nagy, and F.~Allg\"ower.
\newblock {R}eal-time optimization and {N}onlinear {M}odel {P}redictive
  {C}ontrol of {P}rocesses governed by differential-algebraic equations.
\newblock Journal of Process Control, 12(4):577--585, 2002.

\bibitem{Diehl2009}
M.~Diehl and H.J.~Ferreau and N.~Haverbeke.
\newblock {\em Efficient Numerical Methods for Nonlinear MPC and Moving Horizon Estimation},
\newblock In {N}onlinear Model Predictive Control, 384:391--417, 2009.

\bibitem{Feldbaum1961}
A.A.~Feldbaum.
\newblock Dual-control theory (I-IV).
\newblock Automation and Remote Control, 21, pages 1240--1249 and 1453--1464, 1960, as well 22, pages 3--16 and 129--143, 1961.

\bibitem{Ferreau2008}
H.~J. Ferreau, H.~G. Bock, and M.~Diehl.
\newblock {A}n online active set strategy to overcome the limitations of
  explicit {MPC}.
\newblock {\em International Journal of Robust and Nonlinear Control},
  18(8):816--830, 2008.

\bibitem{Filatov2004}
N.M.~Filatov, H.~Unbehauen.
\newblock Adaptive Dual Control. Springer, 2004.

\bibitem{Findeisen2003}
R.~Findeisen, L.~Imsland, F.~Allg\"ower, and B.A.~Foss.
\newblock {\em {S}tate and {O}utput {F}eedback {N}onlinear {M}odel {P}redictive {C}ontrol: {A}n {O}verview},
\newblock European Journal of Control, 9:190--207, 2003.

\bibitem{Griewank2000}
A.~Griewank.
\newblock {\em {E}valuating {D}erivatives, {P}rinciples and {T}echniques of
  {A}lgorithmic {D}ifferentiation}.
\newblock Number~19 in Frontiers in Appl. Math. {SIAM}, Philadelphia, 2000.

\bibitem{Haverbeke2011}
N.~Haverbeke.
\newblock Efficient Numerical Methods for Moving Horizon Estimation.
\newblock Ph.D.~Thesis, KU Leuven, 2011. ISBN~978-94-6018-417-8.

\bibitem{Heirung2012}
T.A.N.~Heirung, B.E.~Ydstie, B.~Foss.
\newblock Towards Dual MPC.
\newblock In Proceedings of the 4th IFAC Nonlinear Model Predictive Control Conference, pp:502--507, Noordwijkerhout (Netherlands), 2012.

\bibitem{Heirung2015}
T.A.N.~Heirung, B.~Foss, B.E.~Ydstie.
\newblock MPC-­‐based dual control with online experiment design.
\newblock Journal of Process Control, 32, 64--76, 2015.

\bibitem{Hernandez2015}
B.~Hernandez, P.~Trodden.
\newblock Persistently exciting tube MPC.
\newblock arXiv:1505.05772v1, 2015.

\bibitem{Hjalmarsson2009}
H.~Hjalmarsson.
System identification of complex and structured systems.
European Journal of Control, 15, 275--310, 2009.

\bibitem{Hovd2005}
M.~Hovd, R.R.~Bitmead.
\newblock Interaction between control and state estimation in nonlinear MPC.
\newblock Modeling, Identification, and Control, 26(3):165--174, 2005.

\bibitem{Houska2011}
B.~Houska, H.J.~Ferreau, and M.~Diehl.
\newblock An auto-generated real-time iteration algorithm for nonlinear {MPC} in the microsecond range.
\newblock Automatica, 47(10):2279-2285, 2011.

\bibitem{Houska2015}
B.~Houska, D.~Telen, F.~Logist, M.~Diehl and J.~Van Impe.
\newblock An Economic Objective for Optimal Experiment Design of Nonlinear Dynamic Processes.
\newblock {\em Automatica}, 51, pp:98--103, 2015.

\bibitem{Houska2016}
B.~Houska, D.~Telen, F.~Logist, J.~Van Impe.
\newblock Self-reflective model predictive control.
\newblock arXiv:1610.03228, 2016.

\bibitem{Larrson2013}
C.A.~Larsson, M.~Annergren, H.~Hjalmarsson, C.R.~Rojas, X.~Bombois, A.~Mesbah, P. Moden.
\newblock Model predictive control with integrated experiment design for output error systems.
\newblock In proceedings of the 2013 European Control Conference (ECC), 3790--3795, July, 2013.

\bibitem{Larsson2015}
C.~Larsson, C.~Rojas, X.~Bombois, H.~Hjalmarsson.
\newblock Experimental evaluation of model predictive control with excitation (MPC-X) on an industrial depropanizer.
\newblock Journal of Process Control, 31:1--16, 2015.

\bibitem{Lee2009}
J.M.~Lee and J.H.~Lee.
\newblock An approximate dynamic programming based approach to dual adaptive control.
\newblock Journal of Process Control, 19(5):859--864, 2009.

\bibitem{Marafioti2014}
G.~Marafioti, R.R.~Bitmead, M.~Hovd.
\newblock Persistently exciting model predictive control.
\newblock International Journal of Adaptive Control and Signal Processing, 28(6):536--552, 2014.

\bibitem{Marchetti2009}
A.~Marchetti, B.~Chachuat, D.~Bonvin.
\newblock Modifier-adaptation methodology for real-time optimization.
\newblock Industrial \& Engineering Chemistry Research, 48(13):6022--6033, 2009.

\bibitem{Mattingley2009}
J.~Mattingley and S.~Boyd.
\newblock Automatic Code Generation for Real-Time Convex Optimization.
\newblock Convex Optimization in Signal Processing and Communications, Y. Eldar and D. Palomar, Eds., Cambridge University Press, 2009.

\bibitem{Mayne2000}
D.Q.~Mayne,  J.B.~Rawlings, C.V.~Rao, and P.Q.~Scokaert.
\newblock Constrained model predictive control: Stability and optimality.
\newblock Automatica, 36(6), 789--814, 2000.

\bibitem{Mayne2006}
D.Q.~Mayne, S.V.~Rakovic, R.~Findeisen, and F.~Allg\"ower.
\newblock {R}obust output feedback model predictive control of constrained
  linear systems.
\newblock {\em Automatica}, 42:1217--1222, 2006.

\bibitem{Mayne2009}
D.Q.~Mayne, S.V.~Rakovic, R.~Findeisen, and F.~Allg\"ower.
\newblock Robust output feedback model predictive control of constrained linear systems:
Time varying case.
\newblock {\em Automatica}, 45:2082--2087, 2009.

\bibitem{Mesbah2015}
A. Mesbah, X. Bombois, M. Forgione, H. Hjalmarsson, P.M.J. Van den Hof.
\newblock Least costly losed-loop performance diagnosis and plant re-identification.
\newblock International Journal of Control, 88(11):2264–2276, 2015.

\bibitem{Nocedal2006}
J.~Nocedal and S.J. Wright.
\newblock {\em {N}umerical {O}ptimization}.
\newblock Springer Series in Operations Research and Financial Engineering.
  Springer, 2 edition, 2006.

\bibitem{Pontryagin1962}
L.S.~Pontryagin, V.G.~Boltyanskii, R.V.~Gamkrelidze, E.F.~Mishchenko.
\newblock The Mathematical Theory of Optimal Processes.
\newblock John Wiley \& Sons, New York, 1962.

\bibitem{Pukelsheim1993}
F.~Pukelsheim.
\newblock Optimal design of Experiments.
\newblock John Wiley \& Sons, Inc., New York, 1993. 

\bibitem{Qin2003}
S.J.~Qin and T.A.~Badgwell.
\newblock A survey of industrial model predictive control technology.
\newblock Control engineering practice, 11(7), 733--764, 2003.

\bibitem{Quirynen2014}
R.~Quirynen, B.~Houska, M.~Vallerio, D.~Telen, F.~Logist, J.~Van Impe, M.~Diehl.
\newblock Symmetric algorithmic differentiation based exact Hessian SQP method and software for economic MPC.
\newblock In Proceedings of the IEEE Conference on Decision and Control (CDC), 2752--2757, 2014.

\bibitem{Rawlings2009}
J.B.~Rawlings and D.Q.~Mayne.
\newblock Model Predictive Control: Theory and Design.
\newblock Madison, WI: Nob Hill Publishing, 2009.

\bibitem{Sager2013}
S.~Sager.
\newblock Sampling decisions in optimum experimental design in the light of Pontryagin's maximum principle.
\newblock SIAM Journal on Control and Optimization 51(4), 3181-3207, 2013.

\bibitem{Shouche2002}
M.S.~Shouche, H.~Genceli, M.~Nikolaou.
\newblock Effect of online optimization techniques on model predictive control and identidification (MPCI).
\newblock Computers \& Chemical Engineering, 26(9):1241--1252, 2002.

\bibitem{Stengel1994}
R.~Stengel.
\newblock Optimal Control and Estimation.
\newblock Dover Publications, New York, 1994.

\bibitem{Telen2013}
D.~Telen, B.~Houska, F.~Logist, E.~Van Derlinden, M.~Diehl, J.~Van Impe.
\newblock Optimal experiment design under process noise using Riccati differential equations.
\newblock Journal of Process Control, 23:613--629, 2013.

\bibitem{Wittenmark1995}
B.~Wittenmark.
\newblock Adaptive dual control methods: An overview.
\newblock In 5th IFAC symposium on Adaptive Systems in Control and Signal Processing Budapest, Hungary, 1995.

\bibitem{Wright1997}
S.J.~Wright.
\newblock Primal-dual interior-point methods.
\newblock SIAM, 1997.

\bibitem{Zacekova2013}
E.~Zacekova, S.~Privara, M.~Pcolka.
\newblock Persistent excitation condition within the dual control framework.
\newblock Journal of Process Control, 23(9):1270--1280, 2013.

\bibitem{Zavala2009}
V.~M. Zavala and L.T. Biegler.
\newblock {T}he {A}dvanced {S}tep {NMPC} {C}ontroller: {O}ptimality,
  {S}tability and {R}obustness.
\newblock {\em Automatica}, 45:86--93, 2009.

\end{thebibliography}

\end{document}